\documentclass[11pt,leqno]{article} 
\usepackage{graphics}
\newtheorem{thm}{Theorem}[section]

\newtheorem{prop}{Proposition}

\newcommand{\beqa}{\begin{eqnarray}}
\newcommand{\eeqa}{\end{eqnarray}}

\newcommand{\pf}{\noindent {\bf Proof:} $\s$ }
\newcommand{\epf}{ \hfill$\diamondsuit$ \medskip}

\newcommand{\beq}{\begin{equation}}
\newcommand{\eeq}{\end{equation}}
\newcommand{\lbl}{\label}
\newcommand{\s}{\; \;}

\newcommand{\la}{\lambda}

\newcommand{\ra}{\rightarrow}
\newcommand{\al}{\alpha}

\newcommand{\p}{\varphi}

\title{Explicit solutions and multiplicity results for some equations with the $p$-Laplacian}

\author{
Philip Korman   \\ 
Department of Mathematical Sciences \\ 
University of Cincinnati \\ 
Cincinnati Ohio 45221-0025 \\
}

\date{}

\begin{document}

\maketitle
\begin{abstract} 
We derive explicit ground state solutions for several equations with the $p$-Laplacian in $R^n$, including (here  $\p (z)=z|z|^{p-2}$, with $p>1$)
\[
\p \left(u'(r)\right)' +\frac{n-1}{r} \p \left(u'(r)\right)+u^M+u^Q=0 \,.
\]
The constant $M>0$ is assumed to be below the critical power, while $Q=\frac{M p-p+1}{p-1}$ is above the critical power. This explicit solution is used to give a multiplicity result, similarly to C.S. Lin and W.-M. Ni \cite{LN}. We also give the $p$-Laplace version of G. Bratu's solution \cite{B}.
\medskip

In another direction, we present a change of variables which removes the non-autonomous term $r^{\al}$ in  
\[
\p \left(u'(r)\right)' +\frac{n-1}{r} \p \left(u'(r)\right)+r^{\al} f(u)=0 \,,
\]
while preserving the form of this equation. In particular, we study singular equations, when $\al <0$. The Coulomb case $\al=-1$ turned out to give the critical power.
 \end{abstract}

\begin{flushleft}
Key words:  Explicit solutions, multiplicity results. 
\end{flushleft}

\begin{flushleft}
AMS subject classification: 35J25, 35J61.
\end{flushleft}

\section{Introduction}
\setcounter{equation}{0}
\setcounter{thm}{0}
\setcounter{lma}{0}

For the equation with the critical exponent (where $u=u(x)$, $ x \in R^n$)
\beq
\lbl{i10}
\Delta u+u^{\frac{n+2}{n-2}}=0
\eeq
there is a well-known explicit solution
\beq
\lbl{i11}
u(x)= \left( \frac{an}{1+\frac{n}{n-2} a^2 r^2} \right)^{\frac{n-2}{2}} \,,
\eeq
see T. Aubin \cite{A} or G. Talenti \cite{T}. Here $r=|x|$, and $a$ is an arbitrary positive  constant. This  explicit solution is very important, for example, it played a central role in the classical paper of H. Brezis and L. Nirenberg \cite{BN}. How does one derive such a solution? Radial solutions of (\ref{i10}) satisfy
\beq
\lbl{i12}
u'' +\frac{n-1}{r}u'+ u^{\frac{n+2}{n-2}}=0 \,, \s  u'(0)=0 \,, \s u'(r)<0 \,.
\eeq
Let us set
\beq
\lbl{i13}
u'=-aru^{\frac{n}{n-2}} \,.
\eeq
Then $u''=-au^{\frac{n}{n-2}}+\frac{n}{n-2} a^2r^2u^{\frac{n+2}{n-2}}$, and using these expressions for $u'$ and $u''$ in (\ref{i12}), we get an algebraic equation for $u$, solving of which leads to the solution in (\ref{i11}). In order for such an approach to work,  the solution $u(r)$ must satisfy the ansatz (\ref{i13}), and it does! 
\medskip

We show that a similar approach produces the explicit solution of C.S. Lin and W.-M. Ni \cite{LN} for the equation
\beq
\lbl{i14}
u'' +\frac{n-1}{r}u'+u^q+u^{2q-1}=0 \,,
\eeq
with $\frac{n}{n-2}<q<\frac{n+2}{n-2}<2q-1$, and some other equations, and for the $p$-Laplace versions of all of these equations. As an application, we state a multiplicity result for  the $p$-Laplace version of (\ref{i14}), similarly to C.S. Lin and W.-M. Ni \cite{LN}.
\medskip

While studying positive solutions of semilinear equations on a ball in $R^n$, we noticed that  for the non-autonomous problem (here $\al>0$, and $a>0$ are constants)
\beq
\lbl{i1}
u'' +\frac{n-1}{r}u'+r^{\al } f(u)=0 \,, \s  u(0)=a \,, \s u'(0)=0 \,,
\eeq
one can prove similar results as for the autonomous case, when $\al =0$. We wondered if the $r^{\al }$  term can be removed by a change of variables.
It turns out that 
the change of variables $t=\frac{r^{1+\al/2}}{1+\al/2}$ transforms the problem (\ref{i1}) into
\beq
\lbl{i2}
u''(t) +\frac{m}{t}u'(t) + f(u(t))=0 \,, \s  u(0)=a \,, \s \frac{du}{dt}(0)=0 \,,
\eeq
with $m=\frac{n-1+\al/2}{1+\al/2}$.
The point here is that this change of variables preserves the Laplacian in the equation. This transformation allows us to get some new multiplicity results for the corresponding Dirichlet problem, including the singular case, when $\al <0$. We present similar results for equations with the $p$-Laplacian.
Such problems, with the $r^{\al }$ term, often arise in applications, for example in modeling of electrostatic micro-electromechanical systems (MEMS), see e.g., J.A. Pelesko \cite{P}, N. Ghoussoub and Y.  Guo \cite{GG},  Z. Guo and J. Wei \cite{GW}.

\section{Some explicit ground state solutions}
\setcounter{equation}{0}
\setcounter{thm}{0}
\setcounter{lma}{0}

For the problem
\beq
\lbl{e1}
u'' +\frac{n-1}{r}u'+f(r,u)=0 \,, \s \mbox{$r>0$}\,, \s u'(0)=0 \,,
\eeq
the crucial role is played by  Pohozhaev's function
\[
P(r)=r^n \left[ {u'}^2(r)+2F(r,u(r)) \right] +(n-2)r^{n-1}u'(r)u(r) \,,
\]
where we denote $F(r,u)=\int _0^u f(r,t) \,dt$. One computes that any solution of (\ref{e1}) satisfies
\beq
\lbl{e2}
\s\s P'(r) = r^{n-1} \left[2nF(r,u(r))-(n-2) u(r)f(r,u(r))+2r F_r(r,u(r)) \right] \,.
\eeq
In case $f(r,u)=u^p$, we have $P'(r)=0$ for $p=\frac{n+2}{n-2}$, $P'(r)<0$ for $p>\frac{n+2}{n-2}$, and $P'(r)>0$ for $p<\frac{n+2}{n-2}$. (Integrating (\ref{e2}), one shows that the Dirichlet problem for (\ref{e1}) on any ball has no solutions if $p>\frac{n+2}{n-2}$.) The  critical exponent $\frac{n+2}{n-2}$  is also the cut-off for the Sobolev embedding. In case $f(r,u)=r^{\al}u^p$, with a constant $\al$,  we have $P'(r)=0$ for $p=\frac{n+2+2\al}{n-2}$, the new critical exponent. Integrating (\ref{e2}), one sees that the Dirichlet problem for the equation (\ref{e3}) below, on any ball, has no solutions if $p>\frac{n+2+2\al}{n-2}$.
\medskip

Let us look for positive ground state solutions of ($n>2$)
\beq
\lbl{e3}
u'' +\frac{n-1}{r}u'+r^{\al}u^{\frac{n+2+2\al}{n-2}}=0 \,, \s \mbox{$r>0$}\,, \s u'(0)=0 \,.
\eeq
Denoting  $p=\frac{n+2+2\al}{n-2}$, we let (observing that $u'(r)<0$)
\beq
\lbl{e4}
u'=-ar^{1+\al} u^{\frac{p+1}{2}}=-ar^{1+\al} u^{\frac{n+\al}{n-2}} \,,
\eeq
where $a>0$ is a constant. Then
\[
u''=-(1+\al)ar^{\al} u^{\frac{p+1}{2}}+ \frac{p+1}{2} a^2 r^{2+2\al} u^{p} \,.
\]
Using these expressions for $u'$ and $u''$ in (\ref{e3}), we get an algebraic expression, which we solve for $u$:
\beq
\lbl{e5}
u(r)=\left[ \frac{a n +a \al}{1+\frac{p+1}{2} a^2 r^{2+\al}} \right] ^{\frac{2}{p-1}}=\left[ \frac{a n +a \al}{1+\frac{n+\al}{n-2} a^2 r^{2+\al}} \right] ^{\frac{n-2}{2+\al}} \,.
\eeq
In order for this function to be a solution of (\ref{e3}), it must satisfy the ansatz (\ref{e4}), which might look  unlikely. But is does, for any constant $a$! By choosing $a$, we can satisfy the initial conditions $u(0)=A$, $u'(0)=0$, for any $A>0$. When $\al =0$, the ground state solution in (\ref{e5}) is the same as the well-known one in (\ref{i11}). 
\medskip

We consider next the problem ($n>2$, $p>1$)
\beq
\lbl{e6}
u'' +\frac{n-1}{r}u'+r^{\al}\left(-u^p+u^{2p-1} \right)=0 \,, \s \mbox{$r>0$}\,, \s u'(0)=0 \,.
\eeq
We set
\beq
\lbl{e7}
u'=-ar^{1+\al} u^p \,,
\eeq
where $a>0$ is a constant.
Then 
\[
u''=-(1+\al)ar^{\al} u^p+a^2pr^{2+2\al}u^{2p-1} \,.
\]
Using these expressions for $u'$ and $u''$ in (\ref{e6}), we obtain
\beq
\lbl{e8}
u(r)=\left[ \frac{a\, n +a\, \al+1}{1+p \, a^2 r^{2+\al}} \right]^{\frac{1}{p-1}} \,.
\eeq
This function satisfies the ansatz (\ref{e7}) provided that 
\beq
\lbl{e9}
a=\frac{p-1}{\alpha -n p+n+2
   p} \,.
\eeq
In order to have $a>0$, we need $p<\frac{n+\al}{n-2}$, and then $2p-1<\frac{n+2+2\al}{n-2}$, i.e., both powers are sub-critical. Conclusion: the function $u(r)$ in (\ref{e8}), with $a$ given by (\ref{e9}) provides a ground state solution for (\ref{e6}).
\medskip

Finally, we  consider  the problem ($n>2$, $p>1$)
\beq
\lbl{e10}
u'' +\frac{n-1}{r}u'+r^{\al}\left(u^p+u^{2p-1} \right)=0 \,, \s \mbox{$r>0$}\,, \s u'(0)=0 \,.
\eeq
Using  the ansatz (\ref{e7}) again, we obtain
\beq
\lbl{e11}
u(r)=\left[ \frac{a\,  n +a\,  \al-1}{1+p \, a^2 r^{2+\al}} \right]^{\frac{1}{p-1}} \,.
\eeq
This function satisfies the ansatz (\ref{e7}) provided that 
\beq
\lbl{e12}
a=\frac{p-1}{np-n-2p-\al} \,.
\eeq
In order to have $a>0$, we need $p>\frac{n+\al}{n-2}$, and then $2p-1>\frac{n+2+2\al}{n-2}$, the  critical exponent. Conclusion: the function $u(r)$ in (\ref{e11}), with $a$ given by (\ref{e12}) provides a ground state solution for (\ref{e10}). In case $\al=0$, this solution was originally found by C.S. Lin and W.-M. Ni \cite{LN}.
\medskip

A similar approach can be tried  for the equations of the form
\beq
\lbl{e14}
u'' +\frac{n-1}{r}u'+A \psi(u)+B \psi(u)\psi'(u)=0 \,, \s \mbox{$r>0$}\,, \s u'(0)=0 \,,
\eeq
where $\psi(u)$ is a given function, with monotone $\psi'(u)$, so that the inverse function ${(\psi')}^{-1}(u)$ exists. Here $A$ and $B$ are given constants. Setting
\beq
\lbl{e15}
u'=-ar \psi(u) \,,
\eeq
with $u''=a^2r^2 \psi(u) \psi'(u)- a\psi(u)$, we obtain from (\ref{e14})
\beq
\lbl{e16}
u(r)={(\psi')}^{-1} \left(\frac{a\, n-A}{a^2r^2+B} \right) \,.
\eeq
This function gives a solution of (\ref{e14}), provided it satisfies (\ref{e15}). If we select here $n=2$, $A=0$, and $\psi(u)=\sqrt2 e^{u/2}$, then the last formula gives
\beq
\lbl{16.1}
u(r)=2 \ln \frac{2 \sqrt2 a}{a^2 r^2+B} \,.
\eeq
One verifies that for any $a>0$, and any $B>0$ the function in (\ref{16.1}) solves 
\[
u''(r)+\frac1r u'(r)+Be^{u(r)}=0 \,, \s u'(0)=0 \,.
\]
This is the famous G. Bratu's \cite{B} solution. It immediately implies the exact count of solutions for the corresponding Dirichlet problem on the unit ball in $R^2$.
\begin{prop}
The problem 
\[
u''(r)+\frac1r u'(r)+Be^{u(r)}=0 \,, \s u'(0)=u(1)=0 
\]
has exactly two solutions for $0<B<2$, exactly one solution for $B=2$, and no solutions if $B>2$.
\end{prop}

\pf
According to the formula (\ref{16.1}), the boundary condition $u(1)=0$ is equivalent to
\[
a^2-2\sqrt2 \, a+B=0 \,.
\]
This quadratic equation has  two solutions for $0<B<2$,  one solution for $B=2$, and none if $B>2$.
\epf

Another  example: the equation 
\[
u'' +\frac{n-1}{r}u'+(n-2)e^u+Be^{2u}=0  \,, \s \mbox{$r>0$}\,, \s u'(0)=0
\]
has a solution $u=\ln \frac{2}{r^2+B}$, for any real $B$.
\medskip

The class of $\psi (u)$, for which this approach works is not  wide. Indeed, writing (\ref{e16}) as $\psi'(u)=\frac{n-A}{r^2+B}$, differentiating this equation, and using (\ref{e15}), we see that $\psi (u)$ must satisfy
\beq
\lbl{e17}
\psi'' (u) \psi (u)=\frac{2}{n-A} {\psi'}^2(u) \,.
\eeq
Solutions of the last equation are exponentials and powers (of $c_1u+c_2$). If $A=0$, a solution of (\ref{e17}) is  $\psi (u)=u^k$, with $k=\frac{n}{n-2}$, which leads to the ground state solution for the critical power $\frac{n+2}{n-2}$, that we considered above.

\section{Explicit ground states in case of the $p$-Laplacian}
\setcounter{equation}{0}
\setcounter{thm}{0}
\setcounter{lma}{0}

For equations with the radial $p$-Laplacian in $R^n$ ($n \geq p$)
\beq
\lbl{l0}
\p \left(u'(r)\right)' +\frac{n-1}{r} \p \left(u'(r)\right)+f(u)=0 \,,
\eeq
Pohozhaev's function
\[
P(r)=r^n \left[ (p-1) \p(  u'(r))u'(r)+pF(u(r)) \right] +(n-p)r^{n-1} \p(  u'(r))u(r) 
\]
was introduced in P. Korman \cite{K}. Here $\p (z)=z|z|^{p-2}$, with $p>1$, and $F(u)=\int_0^u f(t) \, dt$. For the solutions of (\ref{l0}) we have
\[
P'(r)=r^{n-1} \left[npF(u)-(n-p)uf(u)  \right] \,.
\]
Comparing this $P(r)$ to the one in case $p=2$, it was relatively easy for us to make the adjustments, except for the $p-1$ factor, which we found only after a lot of experimentation, using {\em Mathematica}. In case $f(u)=u^q$, one calculates the critical power (when $P'(r)=0$)  to be  $q=\frac{(p-1)n+p}{n-p}$.
\medskip

We look for  positive ground state solutions of ($n > p$)
\beq
\lbl{l1}
\p \left(u'(r)\right)' +\frac{n-1}{r} \p \left(u'(r)\right)+u^q=0 \,, \s u'(0)=0 \,,
\eeq
where $q$ is the critical power $q=\frac{(p-1)n+p}{n-p}$. Then $P'(r)=0$, so that $P(r)=constant=0$, which simplifies as
\beq
\lbl{*+}
r \left[ (p-1) |u'|^p+p\frac{u^{q+1}}{q+1} \right] +(n-p) \p(  u'(r))u(r) =0 \,.
\eeq
By maximum principle, positive solutions of (\ref{l1}) satisfy $u'(r) \leq 0$, for all $r$. In (\ref{*+}) we set ($a>0$ is a constant)
\beq
\lbl{l2}
\p \left(u'(r)\right)=-aru^s(r) \,,
\eeq
with the power $s$ to be specified. Writing (\ref{l2}) as $ \p \left(-u'(r)\right)=aru^s(r)$, or $ \left(-u'(r)\right)^{p-1}=aru^s(r)$, we express $ -u'(r)=a^{\frac{1}{p-1}}r^{\frac{1}{p-1}} u^{\frac{s}{p-1}}(r)$. Then (\ref{*+}) becomes
\beq
\lbl{l3}
(p-1)a^{\frac{p}{p-1}}r^{\frac{p}{p-1}} u^{\frac{sp}{p-1}}+\frac{p}{q+1} u^{q+1} =a(n-p)u^{s+1} \,.
\eeq
We now choose $s$ to get the equal powers of $u$ on the left: $\frac{sp}{p-1}=q+1$, giving 
\[
s=\frac{(q+1)(p-1)}{p}=\frac{n(p-1)}{n-p} \,.
\]
Then solving (\ref{l3}) for $u$, we get
\beq
\lbl{l4}
u(r)=\left[ \frac{a(n-p)}{\frac{n-p}n+(p-1) a^{\frac{p}{p-1}} r^{\frac{p}{p-1}}}          \right]^{\frac{n-p}{p}} \,.
\eeq
One verifies that this $u(r)$ satisfies the ansatz (\ref{l2}) for any $a>0$, and so it gives a ground state solution of (\ref{l1}). (A computation using {\em Mathematica 10} required ``human assistance". {\em Mathematica } calculated $\p \left(u'(r)\right)+aru^s(r)$, and factored the answer, but did not  recognize that one of the factors, $(n-p)^p-p^p
   \left(\frac{n}{p}-1\right)^
   p$, is zero, until it was told that $p>0$.) By choosing $a$, we can satisfy the initial conditions $u(0)=A$, $u'(0)=0$, for any $A>0$.
\medskip

We consider next the equation of Lin-Ni type  with the $p$-Laplacian
\beq
\lbl{l5}
\p \left(u'(r)\right)' +\frac{n-1}{r} \p \left(u'(r)\right)+u^M+u^Q=0 \,.
\eeq
Here $M>p-1$ is a positive constant, and 
\beq
\lbl{l6}
Q=\frac{M p-p+1}{p-1}>M \,.
\eeq
Looking for a positive ground state, we set in (\ref{l5})
\beq
\lbl{l7}
\p \left(u'(r)\right)=-aru^M(r) \,,
\eeq
with the constant $a>0$ to be determined. As above, we express $ -u'(r)=a^{\frac{1}{p-1}}r^{\frac{1}{p-1}} u^{\frac{M}{p-1}}(r)$, so that 
\[
\frac{d}{dr} \p \left(u'(r)\right) =-au^M-arMu^{M-1}u'=-au^M+Ma^{\frac{p}{p-1}}r^{\frac{p}{p-1}} u^Q \,.
\]
Then (\ref{l5}) gives
\beq
\lbl{l8}
u(r)=\left(\frac{a
   n-1}{1+a^{\frac{p}{p-1}} M
   r^{\frac{p}{p-1}}}\right)
   ^{\frac{p-1}{M-p+1}} \,.
\eeq
In order for this function to be a solution of (\ref{l5}), it must  satisfy the ansatz (\ref{l7}). This happens if 
\beq
\lbl{l9}
a=\frac{M-p+1}{M n-p n+n-M p} \,.
\eeq
Observe that $an>1$, provided that both the numerator and denominator are positive in (\ref{l9}), or when 
\beq
\lbl{l10}
M>\frac{np-n}{n-p} \,,
\eeq
which implies  that  $Q> \frac{(p-1)n+p}{n-p}$, the critical power.
Conclusion: the function $u(r)$ in (\ref{l8}), with $a$ from (\ref{l9}), gives a ground state solution of (\ref{l5}),  provided that (\ref{l10}) holds.
\medskip

Similarly to C.S. Lin and W.-M. Ni \cite{LN} the existence of an explicit ground state solution implies a multiplicity result. 

\begin{thm}
Suppose that $p>1$, $n>p$, $M>p-1$, the condition (\ref{l10}) holds, and $Q$ is defined by (\ref{l6}). Then there exists $R_*>0$, so that for $R> R_*$ the problem
\beqa
\lbl{l11}
& \p \left(u'(r)\right)' +\frac{n-1}{r} \p \left(u'(r)\right)+u^M+u^Q=0 \,, \s \mbox{for $0<r<R$} \\ \nonumber
& u'(0)=u(R)=0 \nonumber
\eeqa
has at least two positive solutions.
\end{thm}

\pf
Recall that (\ref{l10}) implies: $p-1<M<\frac{(p-1)n+p}{n-p}<Q$.
Similarly to C.S. Lin and W.-M. Ni \cite{LN}, we employ ``shooting", and consider 
\beqa
\lbl{l12}
& \p \left(u'(r)\right)' +\frac{n-1}{r} \p \left(u'(r)\right)+u^M+u^Q=0 \,, \s \mbox{for $0<r<R$} \\ \nonumber
& u(0)=a \,, \s u'(0)=0  \,. \nonumber
\eeqa
Let $\rho (a)$ denote the first root of $u(r)$, and we say $\rho (a)=\infty$ if $u(r)$ is a ground state solution.
When $a$ is small, one sees by scaling that a multiple of the solution of (\ref{l12}) is an arbitrarily small perturbation of
\beq
\lbl{l12.1}
 \p \left(z'(r)\right)' +\frac{n-1}{r} \p \left(z'(r)\right)+z^M=0 \,, \s z(0)=a \,, \s z'(0)=0 \,.
\eeq
Indeed, setting $u=aw$, and $r=\beta s$, with $\beta =a^{-\frac{M-p+1}{p}}$, the problem  (\ref{l12}) is transformed into 
\[
\frac{d}{ds} \p \left(\frac{dw}{ds}\right) +\frac{n-1}{s} \p \left(\frac{dw}{ds}\right)+w^M+\epsilon w^Q =0 \,, \s w(0)=1 \,, \s w'(0)=0 \,,
\]
with $\epsilon=a^{Q-M}$. Solutions of the last equation are decreasing (while they are positive), and so the $\epsilon w^Q$ term is bounded by $\epsilon w^Q(0)=\epsilon$.
\medskip

For the problem (\ref{l12.1}) it is known (see e.g., \cite{K} or \cite{K2})) that for any $a>0$, the solution $z(r)$ has a unique root, this root tends to infinity as $ a \ra 0$, and $z(r)$ is negative and decreasing  after the root. By the continuity in $\epsilon$, it follows that $\rho (a) < \infty$ for $a$ small, and $\rho (a) \ra \infty$  as $ a \ra 0$. Now denote $A=\{a>0 \s {\big |} \s \rho (a)<\infty \}$. The set $A$ is open, but since we have an explicit ground state, it follows that there exists an interval $(0,\beta) \subseteq A$, with $\beta \notin A$. By the continuous dependence on the initial data, $\lim _{a \uparrow \beta} \rho(a)=\infty$, and the theorem follows, with $R_*=\inf \{ \rho(a) \s {\big |} \s a \in (0,\beta) \}$.
\epf

We now discuss the problem (\ref{l11}) in case $p=2$, when $Q=2M-1$. By scaling, we can transform it to a Dirichlet problem on a unit ball
\beq
\lbl{l13}
\s\s\s\s\s u'' +\frac{n-1}{r}u'+\la \left( u^M+u^{2M-1} \right) =0 \,,\;  0<r<1 \,, \; u'(0)=u(1)=0 \,,
\eeq
with a positive parameter $\la$. The result of C.S. Lin and W.-M. Ni \cite{LN} (extended above), together with the bifurcation theory developed in \cite{KLO}, \cite{OS} and \cite{K2}, implies the existence of a curve of solutions in the $(\la, u(0))$ plane. Along this curve $\la \ra \infty$, when $u(0) \ra 0$, and when $u(0) \ra \beta$. This curve has a {\em horizontal asymptote } at $u(0)=\beta$, see \cite{OS}. Based on the numerical evidence, we conjecture that the solution curve makes exactly one turn to the right in  the $(\la, u(0))$ plane, and it exhausts the set of positive solutions of (\ref{l13}), see Figure $1$. However, the picture changes drastically even if the lower power $M$ is perturbed, see Figure 2. This surprising phenomenon is similar to the one observed by H. Br\'{e}zis and L. Nirenberg \cite{BN}, in case $f(u)=\la u +u^{\frac{n+2}{n-2}}$.

\begin{figure}
\begin{center}
\scalebox{0.8}{\includegraphics{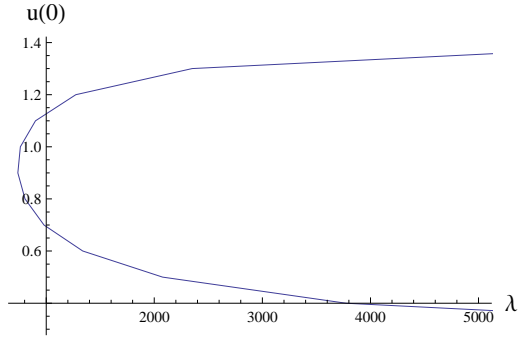}}
\end{center}
\caption{ The solution curve for the problem (\ref{l*})}
\end{figure}
\medskip

\noindent
{\bf Example $1$} We solved numerically the problem (\ref{l13}), with $n=3$, $M=4$, $2M-1=7$
\beq
\lbl{l*}
 u'' +\frac{2}{r}u'+\la \left( u^4+u^7 \right)=0 \,, \s  \; u'(0)=u(1)=0 \,.
\eeq
(See \cite{K2} for the exposition of the shoot-and-scale algorithm that we used.) The solution curve is presented in Figure $1$. Observe that the $\la$'s in this picture are larger than for most other $f(u)$, see \cite{K2}. We have verified this numerical result by an independent computation. Taking an arbitrary point $(\bar \la, \bar u)$ on the solution curve, we solved numerically the initial value problem for the equation in  (\ref{l*}), with $\la =\bar \la$, using the initial conditions $u(0)=\bar u$, $u'(0)=0$. The first root of the solution was always at $r=1$.
\medskip

\noindent
{\bf Example $2$} We solved numerically the problem 
\beq
\lbl{l14}
 u'' +\frac{2}{r}u'+\la \left( u^3+u^7 \right)=0 \,, \s  \; u'(0)=u(1)=0 \,.
\eeq
Compared with the Example $1$, only the lower power is changed from $4$ to $3$. Not only the solution curve, presented in Figure $2$, has a different shape, $\la$'s are now much smaller, while $u(0)$'s go higher. We conjecture that there are still exactly two  positive solutions for $\la$ large enough.

\begin{figure}
\begin{center}
\scalebox{0.8}{\includegraphics{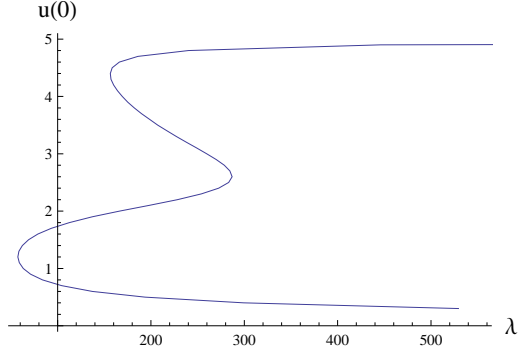}}
\end{center}
\caption{ The solution curve for the problem (\ref{l14})}
\end{figure}
\medskip

We turn next to the $p$-Laplace version of Bratu's equation
\beq
\lbl{l15} 
\p \left(u'(r) \right)'+\frac{n-1}{r} \p \left(u'(r) \right)+Be^u=0 \,,
\eeq
where $\p \left(z \right)=z|z|^{n-1}$ (i.e., $p=n$), and $B>0$ is a constant. Set here 
\[
\p \left(u'(r) \right)=-ar e^{\frac{n-1}{n}u} \,, 
\]
where $a>0$ is a constant. Then $-u'=a ^{\frac{1}{n-1}}r ^{\frac{1}{n-1}}e^{\frac{1}{n}u}$. It follows that
\[
\p \left(u'(r) \right)'=-a e^{\frac{n-1}{n}u}-\frac{n-1}{n}ar e^{\frac{n-1}{n}u}u'=-a e^{\frac{n-1}{n}u}+\frac{n}{n-1} a ^{\frac{n}{n-1}}r ^{\frac{n}{n-1}}e^u \,.
\]
We use these expressions in (\ref{l15}), and solve for $u$:
\beq
\lbl{16.2}
u(r)=n \ln \left( \frac{a \, n}{B+\frac{n}{n-1} a ^{\frac{n}{n-1}}r ^{\frac{n}{n-1}}} \right)\,.
\eeq
One verifies that this function is a solution of (\ref{l15}) for any $a>0$, $B>0$, and $n>1$.
This family of exact solutions  immediately implies the exact count of solutions for the corresponding Dirichlet problem on the unit ball in $R^n$.
\begin{prop}
For the problem 
\[
\p \left(u'(r) \right)'+\frac{n-1}{r} \p \left(u'(r) \right)+Be^u=0 \,, \s u'(0)=u(1)=0 \,,
\]
where $\p \left(z \right)=z|z|^{n-1}$ (i.e., $p=n$), there is a constant $B(n)>0$, so that there are 
 exactly two solutions for $0<B<B(n)$, exactly one solution for $B=B(n)$, and no solutions if $B>B(n)$.
\end{prop}

\pf
According to the formula (\ref{16.2}), the boundary condition $u(1)=0$ is equivalent to $a$ satisfying
\[
\frac{n}{n-1} a ^{\frac{n}{n-1}}+B=n\,  a \,.
\]
On the left we have a convex superlinear function of $a$, so that there is a constant $B=B(n)$, such that  this equation has  two solutions for $0<B<B(n)$,  one solution for $B=B(n)$, and none if $B>B(n)$.
\epf

\section{A change of variables}
\setcounter{equation}{0}
\setcounter{thm}{0}
\setcounter{lma}{0}

For the non-autonomous problem (here $\al$, and $a>0$ are constants)
\beq
\lbl{c1}
u'' +\frac{n-1}{r}u'+r^{\al } f(u)=0 \,, \s  u(0)=a \,, \s u'(0)=0 \,,
\eeq
we present a change of variables which essentially eliminates the non-autonomous term $r^{\al }$ (although it changes the spatial dimension).

\begin{prop}
Let $u(r) \in C^2(0,b) \cap C^1[0,b]$ be a solution of (\ref{c1}), with some $b>0$, and assume that $\al >-1$.
The change of variables $t=\frac{r^{1+\al/2}}{1+\al/2}$ transforms the problem (\ref{c1}) into
\beq
\lbl{c2}
u''(t) +\frac{m}{t}u'(t) + f(u(t))=0 \,, \s  u(0)=a \,, \s \frac{du}{dt}(0)=0 \,,
\eeq
with $m=\frac{n-1+\al/2}{1+\al/2}$.
\end{prop}

\pf
We have $u_r=u_t r^{\al/2}$, $u_{rr}=u_{tt} r^{\al}+\frac{\al}{2} u_t r^{\frac{\al}{2}-1}$, and (\ref{c1}) becomes
\[
u_{tt} r^{\al}+\frac{\al}{2} u_t r^{\frac{\al}{2}-1}+(n-1)u_t r^{\frac{\al}{2}-1}+r^{\al } f(u)=0 \,.
\]
Dividing by $r^{\al}$, we get the equation in (\ref{c2}).
\medskip

To see that $\frac{du}{dt}(0)=0$, we rewrite (\ref{c1}) as $\left(r^{n-1}u' \right)'+r^{\al+n-1}f(u)=0$, and then express
\[
u'(r)=-\frac{1}{r^{n-1}} \int_0^r  z^{\al+n-1}f(u(z)) \, dz \,.
\]
We have
\[
\frac{du}{dt}(0)=\lim _{r \ra 0} \frac{u'(r)}{r^{\al/2}}=-\lim _{r \ra 0} \frac{1}{r^{n-1+\al/2}} \int_0^r  z^{\al+n-1}f(u(z)) \, dz=0 \,.
\]
\epf

Observe that in case $n=2$, we have $m=n-1=1$, which means that the  $r^{\al }$ term is eliminated without changing the dimension. We also remark that for $\al \leq -1$, we do not expect the problem (\ref{c1}) to have solutions of class $ C^2(0,b) \cap C^1[0,b]$, as an explicit example below shows.
\medskip

\noindent
{\bf Example} The problem
\[
u''(t) +\frac{1}{t}u'(t) + e^u=0 \,, \s  u(0)=a \,, \s u'(0)=0
\]
has a solution $u(t)=a-2 \ln \left(1+\frac{e^a
   }{8} t^2\right)$ going back to the paper of G. Bratu \cite{BR}  from 1914, see also J. Bebernes and D. Eberly \cite{B}. (Letting here $a=\ln 8 \left(3 \pm 2 \sqrt2 \right)$, one gets two solutions of the corresponding Dirichlet problem on the unit ball, with $u(1)=0$.) Setting here $t=\frac{r^{1+\al/2}}{1+\al/2}$, we see that 
\beq
\lbl{cc1}
u(r)=a-2 \ln \left(1+\frac{e^a
  }{8
   \left(\frac{\alpha
   }{2}+1\right)^2}  \, r^{\alpha +2} \right)
\eeq is the solution of  the problem
\beq
\lbl{cc2}
u''(r) +\frac{1}{r}u'(r) + r^{\alpha} e^u=0 \,, \s  u(0)=a \,, \s u'(0)=0 \,.
\eeq

This explicit solution is of particular importance for singular equations, when $\al <0$, showing us what to expect for more general nonlinearities than $e^u$. In the  {\em mildly singular } case, when $-1<\al<0$, the function in (\ref{cc1}) is still a solution of (\ref{cc2}), although it is not classical, but only of class $C^{1,1+\al}$. In the  {\em strongly singular } case, when $\al<-1$, the function in (\ref{cc1}) has unbounded derivative as $r \ra 0$. The case of Coulomb potential, when $\al =-1$, is very special. The corresponding solution from (\ref{cc1}) 
\[
 u(r)=a-2 \ln \left( 1+\frac{e^a}{2} r \right)
\]
still satisfies $u(0)=a$, but not $ u'(0)=0$. Instead, we have $ u'(0)=-e^a=-e^{u(0)}$. We see that the initial value problem
\beq
\lbl{cc3}
u''(r) +\frac{1}{r}u'(r) + \frac{1}{r} e^u=0 \,, \s  u(0)=a \,, \s u'(0)=-e^{u(0)}
\eeq
is a natural substitute of the problem (\ref{cc2}) in case of the Coulomb potential. Problems with   the Coulomb potential occur in applications, see J.L. Marzuola et al \cite{m}.
\medskip

We can now extend all of the known multiplicity results for autonomous equations to the non-autonomous equation (\ref{c1}). For example, we have the following result for a cubic nonlinearity, which is based on a similar theorem for $\al =0$ case, see \cite{KLO}, \cite{OS}, \cite{K2}.

\begin{thm}
Assume that $c>2b>0$, and $\al >0$. Then there is a critical $\la _0$, such that for $\la <\la _0$ the problem 
\[
u''+\frac{n-1}{r} u'+\la r^{\al } u(u-b)(c-u)=0, \s r \in (0,1), \s u'(0)=u(1)=0
\]
 has no positive solutions, it has exactly one positive solution at $\la =\la _0$, and there are exactly two positive solutions for $\la >\la _0$. Moreover, all solutions lie on a single smooth solution curve, which for $\la >\la _0$ has two branches, denoted by $u^{-}(r, \la) <u^{+}(r, \la)$, with $u^{+}(r, \la)$ strictly monotone increasing in $\la$, and $\lim _{\la \ra \infty} u^{+}(r, \la)=c$ for all $r \in [0,1)$. For the lower branch, $\lim _{\la \ra \infty} u^{-}(r, \la)=0$ for  $r \ne 0$. (All of the solutions are classical.)
\end{thm}

A similar transformation works for the $p$-Laplace case
\beq
\lbl{c4}
\s\s \p \left(u'(r)\right)' +\frac{n-1}{r} \p \left(u'(r)\right)+r^{\al } f(u(r))=0 \,, \s  u(0)=a \,, \s u'(0)=0 \,,
\eeq
where $\p (z)=z|z|^{p-2}$, with $p>1$.
\begin{prop}
Let $u(r) \in C^2(0,b) \cap C^1[0,b]$ be a solution of (\ref{c4}), with some $b>0$, and assume that $\al >-1$.
The change of variables $t=\frac{r^{1+\al/p}}{1+\al/p}$ transforms the problem (\ref{c4}) into
\beq
\lbl{c5}
\p \left(u'(t)\right)' +\frac{m}{t} \p \left(u'(t)\right) + f(u(t))=0 \,, \s  u(0)=a \,, \s \frac{du}{dt}(0)=0 \,,
\eeq
with $m=\frac{n-1+\al-\al/p}{1+\al/p}$.
\end{prop}

\pf
We have $u_r=u_t r^{\al/p}$, $\p \left(u_r\right)=r^{\al-\al/p}\p \left(u_t\right)$, and 
\[
\frac{d}{dr} \p \left(u_r\right)=(\al-\al/p)r^{\al-\al/p-1}\p \left(u_t\right)+r^{\al-\al/p} \frac{d}{dt}\p \left(u_t\right)  r^{\al/p} \,,
\]
which leads us to (\ref{c5}).
\medskip

To see that $\frac{du}{dt}(0)=0$, we rewrite (\ref{c4}) as $\left(r^{n-1}\p \left(u'\right) \right)'+r^{\al+n-1}f(u)=0$, and then express
\beq
\lbl{c6}
-u'(r)=\left[ \frac{1}{r^{n-1}} \int_0^r  z^{\al+n-1}f(u(z)) \, dz \right]^{\frac{1}{p-1}}\,.
\eeq
We have
\[
-\frac{du}{dt}(0)=\lim _{r \ra 0} \frac{-u'(r)}{r^{\al/p}}=\lim _{r \ra 0} \left[\frac{\left(-u'(r)\right)^{p-1}}{r^{\frac{\al}{p} (p-1)}} \right]^{\frac{1}{p-1}} \,,
\]
and by (\ref{c6})
\[
\lim _{r \ra 0} \frac{\left(-u'(r)\right)^{p-1}}{r^{\frac{\al}{p} (p-1)}} =-\lim _{r \ra 0} \frac{1}{r^{n-1+\al -\al/p}} \int_0^r  z^{\al+n-1}f(u(z)) \, dz=0 \,,
\]
completing the proof.
\epf

In case $n=p$, we have $m=n-1$, which means that the  $r^{\al }$ term is eliminated without changing the dimension.

\end{document}